\documentclass[12pt,a4paper]{article}
\usepackage[latin1]{inputenc}
\usepackage{amsmath}
\usepackage{amsfonts}
\usepackage{amssymb}
\usepackage[left=2.50cm, right=2.50cm, top=2.50cm, bottom=2.50cm]{geometry}
\usepackage[legalpaper,bookmarks=true,colorlinks=true,linkcolor=black,citecolor=blue,urlcolor=blue]{hyperref}
\begin{document}
\begin{center}
\section*{A Functional Principal Component Analysis Approach to Conditional Copula Estimation}
\end{center}
\vspace{3mm}
\vspace{0.5cm}
\begin{center}
 \large{Toihir Soulaimana Djaloud}\\
{\small  \it D\'epartement de Math\'ematique, Universit\'e Alioune Diop, Bambey, S\'en\'egal.}\\
\end{center}
\begin{center}
 \large{Cheikh Tidiane Seck}\footnote{Corresponding author : cheikhtidiane.seck@uadb.edu.sn. BP 30 Bambey Senegal }\\ %\vspace{3mm}
{\small \it D\'epartement de Math\'ematique, Universit\'e Alioune Diop, Bambey, S\'en\'egal.}\\
 \end{center}
\vspace{3mm}
\textbf{Abstract.} The conditional copula model arises when the dependence between random variables is influenced by another covariate. Despite its  importance in modelling  complex dependence structures, there are very few fully nonparametric approaches to  estimate the conditional copula function. In the bivariate setting, the only  nonparametric estimator for the conditional copula  is based on Sklar's Theorem and proposed by Gijbels \textit{et al.} (2011) \cite{gijbels2011conditional}. In this paper, we propose an alternative nonparametric approach %based on functional principal component analysis. We
to construct an estimator for the bivariate conditional copula  from the Karhunen-Lo\`eve representation of a suitably defined conditional copula process.   We establish its consistency and weak convergence  to a limit Gaussian process with explicit covariance function.\\

\noindent\textbf{Keywords.} conditional copula; conditional copula process; Karhunen-Lo\`eve expansion; nonparametric estimation; weak convergence.\\\\

 %In the conditional copula process model, the covariate is considered as a random variable in the conditional dependence structure. %For each fixed value of the covariate, the conditional copula is a deterministic function of its variables, but its overall behavior can be considered random due to the variability of the covariate, hence the notion of a conditional copula process. 
%This paper proposes a nonparametric approach based on Functional Principal Component Analysis (FPCA) to define estimators of the conditional copula using the Karhunen-Lo\`eve representation. The consistency and weak convergence of the estimator will be provided.\\\\
\subsection*{1. Introduction}
Given \((Y_1,Y_2)\)  a pair of random variables and \(X\) a covariate, all defined in the same probability space \((\Omega, \mathcal{A}, \mathbb{P})\). For any fixed value \(x\) of \(X\), the conditional joint distribution of the pair \((Y_1, Y_2)\) given \(X = x\) is defined by:
$$ H_x(y_1,y_2)=\mathbb{P}\left(Y_1 \leq y_1, Y_2 \leq y_2 \mid X = x\right), \qquad\mbox{for all}\quad y_1, y_2 \in \mathbb{R} $$
from which we derive  the conditional marginal distributions of \(Y_1\) and \(Y_2\) given \(X = x\)  as follows: 
\[ F_{1x}(y) = H_x(y, \infty) \qquad \text{and} \qquad F_{2x}(y) = H_x(\infty, y),\qquad\mbox{for any }\quad  y \in \mathbb{R}. \]

If \(F_{1x}\) and \(F_{2x}\) are continuous, then the extension of Sklar's theorem to conditional distributions by \cite{patton2006modelling} ensures the existence and uniqueness of a conditional copula function \(C_x : [0,1]^2 \rightarrow [0,1]\) satisfying:
$$H_x(y_1, y_2) = C_x(F_{1x}(y_1), F_{2x}(y_2)),\qquad\mbox{ for any} \quad (y_1, y_2) \in \mathbb{R}^2$$
or equivalently for all  \((u_1, u_2) \in [0,1]^2\),
$$C_x(u_1, u_2) = H_x\left(F_{1x}^{-1}(u_1), F_{2x}^{-1}(u_2)\right), $$
where \(F_{1x}^{-1}(u) = \inf \{ y \in \mathbb{R} : F_{1x}(y) \geq u \}\) and \(F_{2x}^{-1}(u) = \inf \{ y \in \mathbb{R} : F_{2x}(y) \geq u \}\) denote the conditional quantile functions associated with \(F_{1x}\) and \(F_{2x}\), respectively. Therefore, \(C_x\) characterizes the dependence structure of the pair \((Y_1, Y_2)\) for any fixed value \(x\) of the covariate \(X\).\\

The nonparametric estimation of the conditional copula \(C_x\) from i.i.d. (independent and identically distributed) observations was introduced by \cite{gijbels2011conditional}. Let \(\{(Y_{1i}, Y_{2i}, X_i)\}_{i=1}^{n}\) be an i.i.d. sample from the random vector \((Y_1, Y_2, X)\) taking values in \(\mathbb{R}^3\). An estimator of \(H_x\) is given by
\[ H_{xh}(y_1, y_2) = \sum_{i=1}^{n} w_{ni}(x, h) \mathbb{I}(Y_{1i} \leq y_1, Y_{2i} \leq y_2), \]
where \(w_{ni}(x, h)\) is a weight, \(h\) is a smoothing parameter, and \(\mathbb{I}(\cdot)\) denotes the indicator function. From this, estimators for the conditional marginal distributions \(F_{1x}\) and \(F_{2x}\) can be derived as follows:
\[ F_{1xh}(y) = H_{xh}(y, \infty) \qquad \text{and} \qquad F_{2xh}(y) = H_{xh}(\infty, y). \]

Relying on Sklar's theorem, \cite{gijbels2011conditional} proposed the following estimator:
\[
C_{xh}(u_1, u_2) = H_{xh}\left({F}_{1xh}^{-1}(u_1), {F}_{2xh}^{-1}(u_2)\right),
\]
which can be also expressed as:
\[
C_{xh}(u_1, u_2) = \sum_{i=1}^{n} w_{ni}(x, h) \mathbb{I}(Y_{1i} \leq F_{1xh}^{-1}(u_1), Y_{2i} \leq F_{2xh}^{-1}(u_2)).
\]
However, this estimator is severely biased when at least one of the marginal distributions varies with the covariate \(X\). To overcome this issue, \cite{gijbels2011conditional} applied a transformation on the margins $Y_1$ and $Y_2$ to remove the effect of $X$. This led to the following estimator, which is independent of the marginal distributions $F_1$ and $F_2$:
\[
\tilde{C}_{xh}(u_1, u_2) = \tilde{G}_{xh}\left(\tilde{G}_{1xh}^{-1}(u_1), \tilde{G}_{2xh}^{-1}(u_2)\right),
\]
where \(\tilde{G}_{xh}(u_1, u_2) = \sum_{i=1}^{n} w_{ni}(x, h) \mathbb{I}(\tilde{U}_{1i} \leq u_1, \tilde{U}_{2i} \leq u_2)\) is the empirical  conditional joint distribution based on the transformed observations \((\tilde{U}_{1i}, \tilde{U}_{2i}) = (F_{1X_i g_1}(Y_{1i}), F_{2X_i g_2}(Y_{2i})), i = 1, \dots, n\), here \(F_{jX_i g_j}(Y_{ji})\), for \(j=1,2\), is an estimator of the conditional marginal distribution \(F_{jX_i}\), based on a bandwidth \(g_j = g_{jn} \longrightarrow 0\) as \(n \longrightarrow \infty\). \(\tilde{G}_{1xh}\) and \(\tilde{G}_{2xh}\) are the empirical margins of \(\tilde{G}_{xh}\).

Subsequently, \cite{veraverbeke2011estimation} studied the asymptotic properties of these two estimators by establishing their asymptotic bias and variance, as well as the weak convergence of the corresponding conditional copula processes in the space \(\ell^{\infty}([0,1]^2)\) of bounded real functions on \([0,1]^2\):
$$
\mathbb{C}_{xn}^{\textup{E}}(u_1, u_2) = \sqrt{nh} (C_{xh}(u_1, u_2) - C_x(u_1, u_2)); \quad (0 \leq u_1, u_2 \leq 1)
$$
and
$$
\tilde{\mathbb{C}}_{xn}^{\textup{E}}(u_1, u_2) = \sqrt{nh} (\tilde{C}_{xh}(u_1, u_2) - C_x(u_1, u_2)); \quad (0 \leq u_1, u_2 \leq 1).
$$

\cite{LemyreTaoufikQuessy2015} also obtained the weak convergence of these two processes by using a different approach based on the functional delta method, and extended these results to the case of dependent data.\\
%Let \(\mathbb{D}\) be the space of bivariate distribution functions on \(\mathbb{R}^2\), and \(\Psi:\mathbb{D}\rightarrow \ell^{\infty}([0,1]^2)\) the mapping defined by: \(\Psi(H_x) = H_x \circ (F_{1x}^{-1}, F_{2x}^{-1})\). If \(F_{1x}\) and \(F_{2x}\) are continuous, then the unique copula \(C_x\) associated with \(H_x\) can be written as \(C_x = \Psi(H_x)\), and the process \(\mathbb{C}_{xn}^{\textup{E}}\) can be rewritten as:
%$$
%\mathbb{C}_{xn}^{\textup{E}} = \sqrt{nh} (\Psi(H_{xh}) - \Psi(H_x)).
%$$
%Under certain regularity conditions, \citet{LemyreTaoufikQuessy2015} showed that the process \(\mathbb{H}_{xh} := \sqrt{nh}(H_{xh} - H_x)\) converges weakly to a Gaussian process \(\mathbb{H}_{x}\) with covariance function
%\[
%\text{Cov}\{\mathbb{H}_x(y_1, y_2), \mathbb{H}_x(y'_1, y'_2)\} = K \Gamma_{H_x}(y_1, y_2, y'_1, y'_2),
%\]
%where \(\Gamma_{H_x}(y_1, y_2, y'_1, y'_2) = H_x(y_1 \wedge y'_1, y_2 \wedge y'_2) - H_x(y_1, y_2) H_x(y'_1, y'_2)\) and \(K\) is a positive constant. A direct application of the functional delta method then allows deducing from this result the weak convergence of the processes \(\mathbb{C}_{xn}^{\textup{E}}\) and \(\tilde{\mathbb{C}}_{xn}^{\textup{E}}\) to Gaussian processes.\\
%Moreover, \citet{LemyreTaoufikQuessy2015} extended these results to the case of dependent data.\\

In this paper, we propose a novel nonparametric approach  based on functional principal component analysis to estimate the bivariate conditional copula function. Indeed for any fixed value \(x\) of the covariate $X$,  the conditional copula fucntion \(C_x\) can be viewed as a trajectory of a square-integrable stochastic process \(\{C_X(u,v): 0\le u,v\le 1\}\). The Karhunen-Lo\`eve expansion then allows  to express this trajectory  in terms of the mean function of the process \(C_X\) (which coincides with the partial copula denoted here as \(\bar{C}\))  plus an infinite series of principal component functions.\\

The  paper is organized as follows: In Section 2, we define the conditional copula process $C_X$ on the same probability space as the covariate $X$ and give  its characteristics: mean function, covariance function and covariance operator. In Section 3, we present our novel approach of estimation of the conditional copula $C_x$, for any fixed $x$, which is based on the Karhunen-Lo\`eve representation of the latter. Finally, in Section 4 we establish the consistency and weak convergence of the proposed estimator.

\subsection*{2. Conditional copula process}
Let \((Y_1, Y_2, X)\) be a triplet of random variables defined on a probability space \((\Omega, \mathcal{A}, \mathbb{P})\). We introduce the following random functions for all \(y_1, y_2 \in \mathbb{R}\):
\[ H_X(y_1, y_2) = \mathbb{P}(Y_1 \leq y_1, Y_2 \leq y_2 \mid X), \]
\[ F_{1X}(y_1) = \mathbb{P}(Y_1 \leq y_1 \mid X) \quad \text{and} \quad F_{2X}(y_2) = \mathbb{P}(Y_2 \leq y_2 \mid X). \]
The conditional copula process is the mapping:
\begin{eqnarray}
C_X : \Omega \times [0,1]^2 &\longrightarrow & [0,1] \nonumber \\
(\omega, u, v) &\longmapsto & C_X(\omega, u, v) =: C_{X(\omega)}(u, v)
\end{eqnarray}
such that \(C_{X(\omega)}(u, v) = H_{X(\omega)}(F_{1X(\omega)}^{-1}(u), F_{2X(\omega)}^{-1}(v))\), where \(F_{1X}^{-1}\) and \(F_{2X}^{-1}\) are the generalized inverses of the random functions \(F_{1X}\) and \(F_{2X}\).

For all fixed \(u, v \in [0,1]\), the mapping \(\omega \longmapsto C_{X(\omega)}(u, v)\) is a random variable with a finite second moment, i.e.,
\[ \mathbb{E}[C_X^2(u, v)] = \int_{\Omega} C_{X(\omega)}^2(u, v) \, d\mathbb{P}(\omega) < \infty. \]
For fixed \(\omega \in \Omega\), setting \(X(\omega) = x\), the mapping \((u, v) \longmapsto C_{X(\omega)}(u, v) = C_{x}(u, v)\) is simply the conditional copula associated with the pair \((Y_1, Y_2)\) given \(X = x\). It represents a trajectory or realization of the process \(C_X\); moreover, \(C_x\) is square-integrable on \([0,1]^2\), since
\[ \|C_x\|_2 = \left( \iint_{[0,1]^2} C_{x}^2(u, v) \, du \, dv \right)^{1/2}< \infty. \]
It follows that the process \(C_X\) is square-integrable (i.e., \(\mathbb{E}\|C_X\|_2 < \infty ,\mathbb{E}\) designs the expectation operator). Consequently, \(C_X\) admits a Karhunen-Lo\`eve expansion, i.e, for all \((u, v) \in [0,1]^2\), we have
\begin{equation}\label{KL}
C_X(u, v) = \mathbb{E}C_X(u, v) + \sum_{k=1}^{\infty} \xi_k \varphi_k(u, v),
\end{equation}
where \(\mathbb{E}C_X(u, v)\) is the mean of the process; the functions \(\varphi_k\) are called functional principal components and form an orthonormal basis of \(L^2([0,1]^2)\); the coefficients \(\xi_k\) are defined by \(\xi_k = \langle C_X - \mathbb{E}C_X, \varphi_k \rangle\); they represent random scores associated with the different principal components \(\varphi_k\) and satisfy the following conditions:
\[ \mathbb{E}\xi_k = 0, \quad \mathbb{E}\xi_k^2 = \lambda_k, \quad \text{Cov}(\xi_k, \xi_l) = 0 \; \text{for} \; k \neq l, \]
where \(\lambda_k\) is an eigenvalue of the covariance operator of the process $ C_X $ and \(\varphi_k\) the corresponding eigenfunction. 

The following lemma is useful to compute the mean function of the process $C_X$.
\textbf{Lemma 2.1.}
 Let \(Y_{\jmath}, X\), \(\jmath=1,2\), be random variables with respective distribution functions \(F_{\jmath}\) and \(F\). Then \(F_{\jmath X}(Y_{\jmath}) \sim \mathcal{U}(0,1)\). Then, for all \(t \in [0,1]\), we have \(F_{\jmath X}(Y_{\jmath})\) are uniforms. Furthmore, for all \(x\in \mathbb{R}\), conditionning to \([X=x]\), \(F_{\jmath x}(Y_{\jmath})\) are uniforms.\\
 \textbf{Proof.} For all \(t\in [0,1] \) we have
\begin{eqnarray}\label{cond}
\mathbb{P}(F_{\jmath X}(Y_{\jmath}) \leq t) &=& \int \mathbb{P}(F_{\jmath X}(Y_{\jmath}) \leq t \mid X=x) f_X(x) \, dx\nonumber \\
&=& \int \mathbb{P}(Y_{\jmath} \leq F^{-1}_{\jmath X}(t) \mid X=x) f_X(x) \, dx\nonumber \\
&=& \int F_{\jmath x} \circ F^{-1}_{\jmath x}(t) f_X(x) \, dx\nonumber \\
&=& \int t f_X(x) \, dx = t.
\end{eqnarray}
For \(\jmath=1,2, \)
\begin{eqnarray}
\mathbb{P}(F_{\jmath x}(Y_{\jmath})\leq t|X=x) &=&\mathbb{P}(Y_{\jmath}\leq F_{\jmath x}^{-1}(t)|X=x)\nonumber\\
&=& F_{\jmath x}\circ F_{\jmath x}^{-1}(t)=t.\nonumber\qquad\qquad\qquad\Box
\end{eqnarray}
Using this lemma, we can  express the mean function of the process $ C_X $ as follows:
 \begin{eqnarray*}
\mathbb{E}[C_X(u,v)] &=& \mathbb{E}[H_X(F_{1X}^{-1}(u), F_{2X}^{-1}(v))]\\
&=& \int_{\mathbb{R}} H_x(F_{1x}^{-1}(u), F_{2x}^{-1}(v)) f_X(x) \, dx\\
&=&\int_{\mathbb{R}} \mathbb{P}(F_{1x}(Y_1) \leq u, F_{2x}(Y_2) \leq v \mid X = x) f_X(x) \, dx\\
&=& \mathbb{E}[\mathbb{I}(F_{1x}(Y_1) \leq u, F_{2x}(Y_2) \leq v) \mid X = x],\\
&=& \mathbb{E}[\mathbb{E}[\mathbb{I}(F_{1X}(Y_1) \leq u, F_{2X}(Y_2) \leq v) \mid X]]\\
&=& \mathbb{E}[\mathbb{I}(F_{1X}(Y_1) \leq u, F_{2X}(Y_2) \leq v)]\\
&=&\mathbb{P}(F_{1X}(Y_1) \leq u, F_{2X}(Y_2) \leq v)\\
&=&\bar{C}(u,v).
  \end{eqnarray*}
That is the mean function of the process $C_X$ coincides with the partial copula \( \bar{C} \) representing the joint distribution of the pair  \( (F_{1X}(Y_1), F_{2X}(Y_2)) \).

The covariance function of the process $C_X$ is defined for all  \(\mathbf{u} = (u, v)\) and  \(\mathbf{v} = (u', v')\)  by
\begin{equation}\label{fconv}
\gamma(\mathbf{u}, \mathbf{v}) = \mathbb{E}[(C_X(\mathbf{u}) - \bar{C}(\mathbf{u}))(C_X(\mathbf{v}) - \bar{C}(\mathbf{v}))].
\end{equation}
$\gamma$ represents the kernel of the covariance operator $\Gamma$ defined below  and satisfies the following equation:
\[ \iint_{[0,1]^2} \gamma(\mathbf{u}, \mathbf{v}) \varphi_k(\mathbf{u}) \, d\mathbf{u} = \lambda_k \varphi_k(\mathbf{v}), \]
where   \(\lambda_k\) and  \(\varphi_k\)  are respectively eigenvalue and eigenfunction of the operator $\Gamma$.
In view of  \eqref{KL}, we can  write
\begin{eqnarray}\label{fncn}
\gamma (\mathbf{u},\mathbf{v})=\sum_{k=1}^{\infty}\lambda_k\varphi_k(\mathbf{u})\varphi_k(\mathbf{v}).
\end{eqnarray}
The covariance operator of the process \(C_X\) is defined for all \(f \in L^2([0,1]^2)\) by
$$ \Gamma (f)=\mathbb{E}[\langle C_X-\bar{C},f \rangle (C_X-\bar{C})], $$
or, equivalently in terms of  tensor product,
\[ \Gamma= \mathbb{E}[(C_X-\bar{C}) \otimes (C_X-\bar{C})], \] where for all $f,g,h \in L^2([0,1]^2)$, $(f\otimes g) (h)=<g,h>f $ .

\subsection*{3. A novel approach to conditional copula estimation }
In this section , we propose a new nonparametric estimation  approach using the functional principal components of the conditional copula process $C_X$. 
By conditioning on \(X = x\),  representation \eqref{KL} becomes:
\begin{equation}\label{3.2.3}
C_x(u,v) = \bar{C}(u,v) + \sum_{k=1}^{\infty} \alpha_k(x) \varphi_k(u,v), \qquad (u,v) \in [0,1]^2,
\end{equation}
where \(\alpha_k(x) = \mathbb{E}(\xi_k | X = x)\).

This means that for any fixed value \(x\) of the covariate, the conditional copula \(C_x\) may be decomposed into two terms:  a partial copula \(\bar{C}\) and an additional summation term representing the part of the variation in the dependence structure due to \(x\).

To estimate \(C_x(u,v)\), we have to choose the number $K $ (positive integer) of non-negligible principal components  so that \(C_x\) can be written as:
\begin{equation}\label{sd0}
C_x(u,v) = \bar{C}(u,v) + \sum_{k=1}^{K} \alpha_k(x) \varphi_k(u,v), \qquad (u,v) \in [0,1]^2.
\end{equation}
The choice of \(K\) can be done through several methods:  Scree plot, cummulative varaince percentage (CVP), AIC, BIC and cross-validation methods.
Thus, estimating \(C_x\) involves estimating the partial copula \(\bar{C}\), the functional principal components \(\varphi_k\), and the score coefficients \(\alpha_k(x)\).

Set \(\varepsilon_1 = F_{1X}(Y_1)\), \(\varepsilon_2 = F_{2X}(Y_2)\), and  let \(F_{\varepsilon}\) be the joint distribution of the pair \((\varepsilon_1, \varepsilon_2)\); that is for $(u,v) \in [0,1]^2$,
 \[F_{\varepsilon} (u,v)=\mathbb{P}(\varepsilon_1\leq u,\varepsilon_2\leq v)=\mathbb{P}(F_{1X}(Y_1) \leq u,F_{2X}(Y_2) \leq v).\] 
 We denote the marginal distributions of \(\varepsilon_1\), \(\varepsilon_2\)  by \(F_{\varepsilon_1}\), \(F_{\varepsilon_2}\) respectively. We have
\[ F_{\varepsilon_1}(u) = \mathbb{P}(\varepsilon_1 \leq u) = \mathbb{P}(F_{1X}(Y_1) \leq u )\]
and
\[ F_{\varepsilon_2}(v) = \mathbb{P}(\varepsilon_2 \leq v) = \mathbb{P}(F_{2X}(Y_2) \leq v). \]

According to Lemma 2.1, we have  \( F_{\varepsilon_1}(u)=u\), \(F_{\varepsilon_2}(v)=v\). By Sklar's theorem, there exists a copula function $C_{\varepsilon}$  such that:
$$ F_{\varepsilon}(u,v) = C_{\varepsilon}(F_{\varepsilon_1}(u), F_{\varepsilon_2}(v)).$$
One can observe that $C_{\varepsilon}$ coincides with the partial copula $\bar{C}$. In fact, we have
$$C_{\varepsilon}(u,v)= C_{\varepsilon}(F_{\varepsilon_1}(u), F_{\varepsilon_2}(v)) =F_{\varepsilon}(u,v) =\mathbb{P}(F_{1X}(Y_1) \leq u,F_{2X}(Y_2) \leq v)=\bar{C}(u,v).$$
%=\mathbb{P}(\varepsilon_1\leq u,\varepsilon_2\leq v)
Thus, an estimator of the partial copula $\bar{C}$ can be deduced from that of the distribution $F_{\varepsilon}$. Consider an i.i.d. sample \((Y_{11}, Y_{21}, X_1), \dots, (Y_{1n}, Y_{2n}, X_n)\) of the random vector \((Y_1, Y_2, X)\). For any fixed value \(x\) of the covariate \(X\), if the conditional marginal distributions \(F_{1x}\) and \(F_{2x}\) are known, then \( \varepsilon_{1i} = F_{1X_i}(Y_{1i})\) and \(\varepsilon_{2i} = F_{2X_i}(Y_{2i})\), for \(i = 1, \dots, n\), are true  observations associated with the pair $(\varepsilon_{1} , \varepsilon_{2} )$, with  empirical distribution function  
\[ F_{\varepsilon, n}(u, v) = \frac{1}{n} \sum_{i=1}^n \mathbb{I}\{\varepsilon_{1i} \le u, \varepsilon_{2i} \le v\}. \]
Then an estimator of \(\bar{C}\) may be defined as
\begin{equation}\label{sd0} \bar{C}_n(u, v) = F_{\varepsilon, n}(u, v). \end{equation}

If for any fixed value \(x\), the marginal distributions \(F_{1x}\) and \(F_{2x}\) are unknown, we  deal with \(\tilde{\varepsilon}_1 = \hat{F}_{1X}(Y_1)\), \(\tilde{\varepsilon}_2 = \hat{F}_{2X}(Y_2)\), where \(\hat{F}_{1X}\) and \(\hat{F}_{2X}\) are consistent estimators of the margins \(F_{1X}\) and \(F_{2X}\).  Let \(F_{\tilde{\varepsilon}}\) be the joint distribution of the pair \((\tilde{\varepsilon}_1, \tilde{\varepsilon}_2)\), with margins: % of \(\tilde{\varepsilon}_1\) and \(\tilde{\varepsilon}_2\) are given by:
\[ F_{\tilde{\varepsilon}_1}(y_1) = \mathbb{P}(\tilde{\varepsilon}_1 \leq y_1) = \mathbb{P}(\hat{F}_{1X}(Y_1) \leq y_1),\quad y_1\in\mathbb{R} \]
and 
\[ F_{\tilde{\varepsilon}_2}(y_2) = \mathbb{P}(\tilde{\varepsilon}_2 \leq y_2) = \mathbb{P}(\hat{F}_{2X}(Y_2) \leq y_2), \quad y_2\in\mathbb{R}\]
 which are assumed to be continuous. Then, by using Sklar's theorem, there exists a unique partial copula \(\tilde{\bar{C}}\) such that for all $y_1,y_2\in\mathbb{R}$
\[ F_{\tilde{\varepsilon}}(y_1,y_2) = \tilde{\bar{C}}(F_{\tilde{\varepsilon}_1}(y_1), F_{\tilde{\varepsilon}_2}(y_2)) \]
or equivalently
\[ \tilde{\bar{C}}(u, v) = F_{\tilde{\varepsilon}, n}\big(F_{\tilde{\varepsilon}_1}^{-1}(u), F_{\tilde{\varepsilon}_2}^{-1}(v) \big),\quad (u,v)\in[0,1]^2.\]
Thus, an estimator of \(\tilde{\bar{C}}\) can be defined as:
\begin{equation} \label{ss0}\tilde{\bar{C}}_n(u, v) = F_{\tilde{\varepsilon}, n}\big(F_{\tilde{\varepsilon}_1, n}^{-1}(u), F_{\tilde{\varepsilon}_2, n}^{-1}(v) \big), \end{equation}
where \(F_{\tilde{\varepsilon}, n}\) is the empirical distribution function associated with the sequence of the pseudo-observations \((\tilde{\varepsilon}_{1i}=\hat{F}_{1X_i}(Y_1), \tilde{\varepsilon}_{2i}=\hat{F}_{2X_i}(Y_2))\), $i=1,\cdots,n$, i.e.,
\[ F_{\tilde{\varepsilon}, n}(y_1,y_2) = \frac{1}{n} \sum_{i=1}^n \mathbb{I}\{\tilde{\varepsilon}_{1i} \le y_1, \tilde{\varepsilon}_{2i} \le y_2\}, \]
and \(F_{\tilde{\varepsilon}_1, n}\) and \(F_{\tilde{\varepsilon}_2, n}\) denote the corresponding empirical marginals.

To estimate the functional principal components \(\varphi_k\), we need an estimator of the covariance operator \(\Gamma\). For this, consider an i.i.d. sample \(C_{X_1}, \ldots, C_{X_n}\)  of random functions from the process \(C_X\). An estimator of the covariance function \(\gamma\) may be defined as 
\[ \hat{\gamma}(\mathbf{u}, \mathbf{v}) = \frac{1}{n} \sum_{i=1}^n (C_{X_i}(\mathbf{u}) - \tilde{\bar{C}}_n(\mathbf{u}))(C_{X_i}(\mathbf{v}) - \tilde{\bar{C}}_n(\mathbf{v})), \]
for all \(\mathbf{u} = (u, v)\) and \(\mathbf{v} = (u' ,v')\) in \([0,1]^2\).  
Denote respectively \(\hat{\lambda}_k\) and \(\hat{\varphi}_k\)  the estimators of the eigenvalue \(\lambda_k\)  and eigenfunction \(\varphi_k\). Then we have 
\begin{equation}\label{op}
\hat{\Gamma} (\hat{\varphi}_k)(\mathbf{u})=\iint_{[0,1]^2} \hat{\gamma}(\mathbf{u}, \mathbf{v}) \hat{\varphi}_k(\mathbf{v}) \, \mathrm{d}\mathbf{v} = \hat{\lambda}_k \hat{\varphi}_k(\mathbf{u}),  \end{equation}
 where \( \hat{\Gamma} \) is an estimator of the operator $\Gamma$.

 Under some mild regularity conditions, we have the following statements (see, e.g., Lemma 4.3 in \cite{liu2016functional}), which guarantee the consistency of $\hat{\lambda}_k $ and $\hat{\varphi}_k $:  for all $k = 1, 2, \dots, K$,
\begin{eqnarray}
|\hat{\lambda}_k - \lambda_k| &=& O_P(n^{-1/2}), \\
\Vert\hat{\varphi}_k - \varphi_k\Vert_2 &=& O_P(n^{-1/2}), \label{sd01} \\
\sup_{u, v \in [0, 1]^2}|\hat{\varphi}_k(u, v) - \varphi_k(u, v)| &=& O_P(n^{-1/2}). \label{sd1}
\end{eqnarray}
 %of the covariance operator \(\Gamma\).

The score coefficient $\alpha_k(x) = \mathbb{E}(\xi_k | X=x)$ given by \eqref{sd0} can be estimated through nonparametric regression method, i.e.,
\begin{equation}
 \hat{\alpha}_k(x) = \sum_{i=1}^n w_{ni}(x)\hat{\xi}_{ki},
\end{equation}
with $\hat{\xi}_{ki} = \left< C_{X_i} - \tilde{\bar{C}}_n, \hat{\varphi}_k \right>$ and $w_{ni}(\cdot)$ being a weight function. Assuming that the function $\alpha_k$ is twice continuously differentiable and that the density function $f_X$ of the covariable $X$ is continuous at $x$, we have the following convergence in probability:
\begin{equation}\label{sd2}
 |\hat{\alpha}_k(x) - \alpha_k(x)| = O_P(\tau_n), \qquad \text{where} \quad \tau_n \longrightarrow 0, \quad n \longrightarrow \infty.
\end{equation}
Finally, we define our estimator of the conditional copula function $C_x$ as 
\begin{equation}\label{sd5}
\hat{C}_x(u, v) = \tilde{\bar{C}}_n(u, v) + \sum_{k=1}^{K}\hat{\alpha}_k(x)\hat{\varphi}_k(u, v), \qquad (u, v) \in [0, 1]^2.
\end{equation}

\subsection*{4. Consistency and weak convergence of the estimator}
We begin by proving the consistency of the estimator $\hat{C}_x$. The following proposition gives the uniform convergence in probability. \\\\
\textbf{Proposition 4.1.}
Assume that $\alpha_k(x)$ is twice continuously differentiable and the covariable density $f_X$ is continuous at every $x$. Then for all $(u,v) \in [0,1]^2$, we have for $\tau_n$ satisfying relation \eqref{sd2}, as $n \to \infty$,
\begin{equation}
\sup_{u,v \in [0,1]}|\hat{C}_x(u,v) - C_x(u,v)| = O_P(\tau_n \vee n^{-1/2}),
\end{equation}
where $a \vee b = \max(a,b)$.

\textbf{Proof.}
The difference $\hat{C}_x - C_x$ can be decomposed as follows: for all $(u,v) \in [0,1]^2,$
\begin{eqnarray}
\hat{C}_x(u,v) - C_x(u,v) &=& \tilde{\bar{C}}_n(u,v) - \bar{C}(u,v) \nonumber \\
&+& \sum_{k=1}^K (\hat{\alpha}_k(x) - \alpha_k(x))(\hat{\varphi}_k(u,v) - \varphi_k(u,v)) \nonumber \\
&+& \sum_{k=1}^K (\hat{\alpha}_k(x) - \alpha_k(x)) \varphi_k(u,v) \\
&+& \sum_{k=1}^K \alpha_k(x) (\hat{\varphi}_k(u,v) - \varphi_k(u,v)) \nonumber.
\end{eqnarray}
Taking the supremum,
\begin{eqnarray}
\sup_{0 \leq u,v \leq 1} |\hat{C}_x(u,v) - C_x(u,v)| &\leq& \sup_{0 \leq u,v \leq 1} |\tilde{\bar{C}}_n(u,v) - \bar{C}(u,v)| \nonumber \\
&+& \sum_{k=1}^K |\hat{\alpha}_k(x) - \alpha_k(x)| ~ \sup_{0 \leq u,v \leq 1} |\hat{\varphi}_k(u,v) - \varphi_k(u,v)| \nonumber \\
&+& \sum_{k=1}^K |\hat{\alpha}_k(x) - \alpha_k(x)| ~ \|\varphi_k\|_2 \nonumber \\
&+& \sum_{k=1}^K |\alpha_k(x)| ~ \sup_{0 \leq u,v \leq 1} |\hat{\varphi}_k(u,v) - \varphi_k(u,v)| \nonumber.
\end{eqnarray}

To study the term $|\tilde{\bar{C}}_n(u,v) - \bar{C}(u,v)|$, we interpose $\bar{C}_n$ (which corresponds to the estimator obtained when the margins are known). We have
$$|\tilde{\bar{C}}_n(u,v) - \bar{C}(u,v)| \leq |\tilde{\bar{C}}_n(u,v) - \bar{C}_n(u,v)| + |\bar{C}_n(u,v) - \bar{C}(u,v)|.$$

According to Proposition 1.8 in \cite{deheuvels2009}, we have
$$\sup_{0 \leq u,v \leq 1} |\tilde{\bar{C}}_n(u,v) - \bar{C}_n(u,v)| = \frac{1}{n},$$
and the weak law of large numbers shows that $\bar{C}_n(u,v) - \bar{C}(u,v) = o_P(1)$, thus
\begin{equation}\label{deH}
\sup_{0 \leq u,v \leq 1} |\tilde{\bar{C}}_n(u,v) - \bar{C}(u,v)| = \frac{1}{n} + o_P(1).
\end{equation}

Combining relations \eqref{sd1}, \eqref{sd2},  \eqref{deH}, and the fact that the coefficients $\alpha_k(x)$ are finite for any fixed value $x$,  we finally obtain 
$$\sup_{0 \leq u,v \leq 1} |\hat{C}_x(u,v) - C_x(u,v)| = O_P(\tau_n \vee n^{-1/2}). \qquad \qquad \Box$$

Now, we investigate the weak convergence of  the process $\{\sqrt{n}[\hat{C}_x(u,v) - C_x(u,v)]: 0 \leq u,v \leq 1\}$ in $L^2([0,1]^2)$. Let's decompose
\begin{eqnarray*}
 \sqrt{n}[\hat{C}_x(u,v) &-& C_x(u,v)] =:  \sqrt{n}[\tilde{\bar{C}}_n(u,v) - \bar{C}(u,v)] + \textup{I}(u,v) + \textup{II}(u,v)  \\
&+& \sum_{k=1}^K \alpha_k(x) \sqrt{n}[\hat{\varphi}_k(u,v) - \varphi_k(u,v)] \\
&=& \sqrt{n}[\tilde{\bar{C}}_n(u,v) - \bar{C}_n(u,v)] + \sqrt{n}[{\bar{C}}_n(u,v) - \bar{C}(u,v)]  \\
&+& \sum_{k=1}^K \alpha_k(x) \sqrt{n}[\hat{\varphi}_k(u,v) - \varphi_k(u,v)] + \textup{I}(u,v) + \textup{II}(u,v),
\end{eqnarray*} 
with 
$$\textup{I}(u,v) = \sqrt{n} \sum_{k=1}^K (\hat{\alpha}_k(x) - \alpha_k(x))(\hat{\varphi}_k(u,v) - \varphi_k(u,v))$$ 
and 
$$\textup{II}(u,v) = \sqrt{n} \sum_{k=1}^K \alpha_k(x)(\hat{\varphi}_k(u,v) - \varphi_k(u,v)).$$
We can show that the terms $\textup{I}(u,v)$ and $\textup{II}(u,v)$ are $O_P(\cdot)$ (big-O in probability). Indeed, considering \eqref{sd1} and \eqref{sd2}, we can write for all $k=1, \ldots, K$:
\begin{eqnarray*}
 |(\hat{\alpha}_k(x) - \alpha_k(x))(\hat{\varphi}_k(u,v) - \varphi_k(u,v))| &\leq& | \hat{\alpha}_k(x) - \alpha_k(x)| \times O_P(n^{-1/2}) \\
 &=& O_P(\tau_n n^{-1/2}).
\end{eqnarray*} 
Since $K < \infty$, we deduce that:
\begin{equation}\label{sd3}
\textup{I}(u,v) = O_P(\tau_n).
\end{equation}
Similarly, for all $k=1, \ldots, K$, we have, with $\|\varphi_k\|_2 = 1$,
$$|(\hat{\alpha}_k(x) - \alpha_k(x)) \varphi_k(u,v)| \leq |\hat{\alpha}_k(x) - \alpha_k(x)| \|\varphi_k\|_2 = O_P(\tau_n).$$
This implies that
\begin{equation}\label{sd4}
\textup{II}(u,v) = O_P(\sqrt{n}\tau_n). 
\end{equation}
In view of  \eqref{sd3} and \eqref{sd4}, we can write 
\begin{eqnarray*}
\sqrt{n}[\hat{C}_x(u,v) - C_x(u,v)] &=& 
\sqrt{n}[\tilde{\bar{C}}_n(u,v) - \bar{C}_n(u,v)] + O_P(\sqrt{n}\tau_n) \\
&+& \sqrt{n}[{\bar{C}}_n(u,v) - \bar{C}(u,v)] + O_P(\tau_n) \\
&+& \sum_{k=1}^K \alpha_k(x) \sqrt{n}[\hat{\varphi}_k(u,v) - \varphi_k(u,v)].
\end{eqnarray*}
By using Proposition 1.8 in \cite{deheuvels2009}, we can infer that 
\begin{equation}\label{gibel}
\sup_{u,v \in [0,1]} \sqrt{n}[\tilde{\bar{C}}_n(u,v) - \bar{C}_n(u,v)] = o_P(1).
\end{equation}
Finally, we get
\begin{eqnarray}
\sqrt{n}[\hat{C}_x(u,v) - C_x(u,v)] &=& \sqrt{n}[\bar{C}_n(u,v) - \bar{C}(u,v)] + o_P(1) + O_P(\tau_n) \nonumber \\
&+& \sum_{k=1}^K \alpha_k(x) \sqrt{n}[\tilde{\varphi}_k(u,v) - \varphi_k(u,v)] + O_P(\sqrt{n}\tau_n) \nonumber \\
&:=& \mathbb{G}_n(u,v) + \mathbb{Z}_{xn}(u,v) + O_P(\sqrt{n}\tau_n), \label{dec2}
\end{eqnarray}
where
$$ \mathbb{G}_n(u,v)=\sqrt{n}[\bar{C}_n(u,v) - \bar{C}(u,v)]\qquad \mbox{and}\qquad\mathbb{Z}_{xn}(u,v)= \sum_{k=1}^K \alpha_k(x) \sqrt{n}[\tilde{\varphi}_k(u,v) - \varphi_k(u,v)]. $$
%In decomposition \eqref{dec2}, the first term on the right-hand side can be represented as follows:
The process $\mathbb{G}_n(u,v)$ coincides with the bivariate uniform empirical process. We can rewrite it as 
\begin{equation}\label{EMP}
\mathbb{G}_n(u,v) = \frac{1}{\sqrt{n}} \sum_{i=1}^n \psi(\varepsilon_{1i}, \varepsilon_{2i}, u, v),\quad \mbox{ for all}\quad (u,v)\in[0,1]^2,
\end{equation}
where $\varepsilon_{1i} = F_{1X_i}(Y_{1i})$, $\varepsilon_{2i} = F_{2X_i}(Y_{2i})$, $i = 1, \dots, n$, and $\psi(x,y,u,v) = \mathbb{I}(x \le u, y \le v) - \bar{C}(u,v)$. $\mathbb{G}_n$  converges in distribution to  a Brownian bridge $\mathbb{G}$ associated to the partial copula $\bar{C}$, with covariance function:
$$ {\rm Cov}(\mathbb{G}(u,v), \mathbb{G}(u',v')) = \bar{C}(u \wedge u', v \wedge v') - \bar{C}(u,v) \bar{C}(u',v').$$

Now let us study the behavior of the second process $\mathbb{Z}_{xn}(u,v)$ in the right-hand side of equality \eqref{dec2}. Recall that the covariance operator $\Gamma$ of the process $C_X$ is estimated by
$$\hat{\Gamma} = \frac{1}{n} \sum_{i=1}^n (C_{X_i} - \tilde{\bar{C}}_n) \otimes (C_{X_i} - \tilde{\bar{C}}_n).$$
$\hat{\Gamma}$ and $\Gamma$ are viewed as elements of $\mathcal{S}$, the space Hilbert-Schmidt operators on $L^2([0,1]^2)$, equipped with the norm $\Vert T \Vert_{\mathcal{S}}=\sqrt{\left<T,T\right>_{\mathcal{S}}},\; T\in \mathcal{S}$,  where for all $T_1, T_2\in\mathcal{S}$,
$$ \left<T_1, T_2\right>_{\mathcal{S}}=\sum_{j\geq 1} \left<T_1 u_j, T_2 u_j\right>$$
for any  complete orthonormal system $\{u_j, j\geq 1\} $ in $L^2([0,1]^2)$ and $\left<\cdot,\cdot\right>$ stands for the inner product in $L^2([0,1]^2)$.
Let $Z_n := \sqrt{n} [\hat{\Gamma} - \Gamma]$. The following proposition gives the weak convergence of the empirical mean $\hat{C} = \frac{1}{n} \sum_{i=1}^n C_{X_i}$ in $L^2([0,1]^2)$ and the weak convergence of the operator $Z_n$ to a centered Gaussian Hilbert-Schmidt operator ${Z}$ in $\mathcal{S}$. It will be used to derive the asymptotic behavior fo the empirical eigenfunctions $\hat{\varphi}_k$. \\

%with covariance $\Lambda = \mathbb{E}(Z \otimes Z)$.\\ Let $\bar{C}$ denote the partial copula.\\
\noindent\textbf{Proposition 4.2.}
Assume that $\{C_{X_i}\}_{i=1}^n$ is an i.i.d. sequence of random functions such that $\mathbb{E}\|C_{X_1}\|^2_2 < \infty$, then
\begin{equation}\label{conv-proc}
\sqrt{n}(\hat{C} - \bar{C}) \stackrel{d}{\longrightarrow} \mathcal{N}(0,\Gamma),
\end{equation} 
where  $\Gamma = \mathbb{E}[(C_{X_1} - \bar{C}) \otimes (C_{X_1} - \bar{C})]$ and $\bar{C}$ is the mean of the process $C_X$. Furthermore, if $\mathbb{E}\|C_{X_1}\|^4_2 < \infty$, then
\begin{equation}\label{zn}
Z_n = \sqrt{n}(\hat{\Gamma} - \Gamma) \stackrel{d}{\longrightarrow} {Z},
\end{equation}
where  $Z = [(C_{X_1} - \bar{C}) \otimes (C_{X_1} - \bar{C}) - \Gamma] \in \mathcal{S}$ is a centered Gaussian operator with covariance $\Lambda = \mathbb{E}(Z \otimes Z)$.\\

\textbf{Proof.}
The asymptotic normality of the empirical mean $\hat{C}$ given by \eqref{conv-proc} follows directly from the central limit theorem in the Hilbert space  $L^2([0,1]^2)$. To establish \eqref{zn}, we decompose
\begin{equation}\label{gamma}
\sqrt{n}(\hat{\Gamma} - \Gamma) = \sqrt{n}(\hat{\Gamma} - \tilde{\Gamma}) + \sqrt{n}(\tilde{\Gamma} - \Gamma),
\end{equation}
where
$$\tilde{\Gamma} = \frac{1}{n} \sum_{i=1}^n (C_{X_i} - \bar{C}) \otimes (C_{X_i} - \bar{C}).$$
By applying the central limit theorem in the space $\mathcal{S}$ of  Hilbert-Schmidt operators (which is an Hilbert space) the second term on the right-hand side of equality \eqref{gamma} converges wealky to a Gaussian Hilbert-Schmidt operator, with covariance
$$\Lambda = \mathbb{E}\left[\left\lbrace (C_{X_1} - \bar{C}) \otimes (C_{X_1} - \bar{C}) - \Gamma \right\rbrace\otimes \left\lbrace (C_{X_1} - \bar{C}) \otimes (C_{X_1} - \bar{C}) - \Gamma \right\rbrace \right].$$
For more details, see Theorem 12.3.1 in \cite{kokoszka2017introduction}.

It remains now to show that the other term $\sqrt{n}(\tilde{\Gamma} - \hat{\Gamma}) \stackrel{P}{\longrightarrow} 0$. To do this, we write
\begin{align*}
\sqrt{n}(\tilde{\Gamma} - \hat{\Gamma}) &= \frac{1}{\sqrt{n}} \sum_{i=1}^n \left[(C_{X_i} - \bar{C}) \otimes (C_{X_i} - \bar{C}) - (C_{X_i} - \tilde{\bar{C}}_n) \otimes (C_{X_i} - \tilde{\bar{C}}_n)\right] \\
&= \frac{1}{\sqrt{n}} \sum_{i=1}^n \left[-\bar{C} \otimes C_{X_i} - C_{X_i} \otimes \bar{C} + \bar{C} \otimes \bar{C} + \tilde{\bar{C}}_n \otimes C_{X_i} + C_{X_i} \otimes \tilde{\bar{C}}_n - \hat{\bar{C}}_n \otimes \tilde{\bar{C}}_n \right] \\
&= \frac{1}{\sqrt{n}} \left[-n \bar{C} \otimes \hat{C} - n \hat{C} \otimes \bar{C} + n \bar{C} \otimes \bar{C} + n \tilde{\bar{C}}_n \otimes \hat{C} + n \hat{C} \otimes \tilde{\bar{C}}_n - n \tilde{\bar{C}}_n \otimes \tilde{\bar{C}}_n \right] \\
&= \frac{1}{\sqrt{n}} \left[\sqrt{n}(\hat{C} - \bar{C}) \otimes \sqrt{n}(\hat{C} - \bar{C}) - \sqrt{n}(\tilde{\bar{C}}_n - \hat{C}) \otimes \sqrt{n}(\tilde{\bar{C}}_n - \hat{C})\right].
\end{align*}
We also have
\begin{eqnarray}
\sqrt{n}(\tilde{\bar{C}}_n - \hat{C}) &=& \sqrt{n}(\tilde{\bar{C}}_n - \bar{C}) + \sqrt{n}(\bar{C} - \hat{C}) \nonumber \\
&=& \sqrt{n}(\tilde{\bar{C}}_n - \bar{C}_n) + \sqrt{n}(\bar{C}_n - \bar{C}) + \sqrt{n}(\bar{C} - \hat{C}) \nonumber 
\end{eqnarray}
According to relations \eqref{gibel} and \eqref{conv-proc}, we can write that $\sqrt{n}(\tilde{\bar{C}}_n - \bar{C}_n) = o_P(1)$, $\sqrt{n}(\hat{C} - \bar{C}) = O_P(1)$, and $\sqrt{n}(\bar{C}_n - \bar{C}) = O_P(1)$. Therefore, $\sqrt{n}(\hat{\bar{C}} - \hat{C}) = O_P(1)$, implying $\sqrt{n}(\tilde{\Gamma} - \hat{\Gamma}) = O_P(1/\sqrt{n}) \stackrel{P}{\longrightarrow} 0$. \qquad\qquad$\Box$\\

Now, for all $k=1,\ldots, K$ and $n\ge 1$, let's introduce the random function
$$ T_{k,n}(u,v) = \sum_{j:j\neq k} (\lambda_k - \lambda_j)^{-1} \left< Z_n, \varphi_k \otimes \varphi_j \right> \varphi_j(u,v) ,\qquad  (u,v)\in [0,1]^2, $$
where $\{\varphi_j\}$ denotes an orthonormal basis of eigenfunctions of the space $L^2([0,1]^2)$. Notice that for each $k$ and $n$ $T_{k,n}\in L^2([0,1]^2)$. Since $Z_n$ converges in distribution to $Z$ by Proposition 4.2,  the continuous mapping theorem implies that  for all fixed $k$ the sequence $(T_{k,n},n\geq 1)$ also converges in distribution to a random function $T_k$ in $L^2([0,1]^2)$ such that
$$ T_k(u,v) = \sum_{j:j\neq k} (\lambda_k - \lambda_j)^{-1} \left< Z, \varphi_k \otimes \varphi_j \right> \varphi_j(u,v). $$

The following lemma allows us to represent the process $\mathbb{Z}_{xn}$ in terms of the processes $T_{k,n}$, $k=1,\cdots, K$. We need the following hypothesis:

(H.1) -- The first $K$ eigenvalues of the covariance operator $\Gamma$ are distinct and ordered as follows:
$$ \lambda_1 > \lambda_2 > \cdots > \lambda_K > \lambda_{K+1} \geq \cdots .$$

\noindent\textbf{Lemma 4.1.} Suppose assumptions of  Proposition 4.2 and hypothesis (H.1) are true, then we have
$$\Vert\sqrt{n}[\hat{\varphi}_k - \varphi_k] - T_{k,n} \Vert_2 \stackrel{P}{\longrightarrow} 0, \quad n \longrightarrow \infty.$$
\textbf{Proof.} (see, Kokoska and Reimherr (2013) )\\
 From Lemma 4.1, we can deduce that for every $(u,v) \in [0,1]^2$,
 \begin{equation} 
\mathbb{Z}_{xn}(u,v)=\sum_{k=1}^K \alpha_k(x)T_{k,n}(u,v)+ o_P(1).
\end{equation}
 But for each $k=1,\ldots,K$,  the process $\{T_{k,n}(u,v): 0\le u,v\le1\}$ converges in distribution to a Gaussian process $\{T_{k}(u,v): 0\le u,v\le1\}$, with covariance function given, for all $\textbf{u}=(u,v)$, $\textbf{v}=(u',v')$ in $[0,1]^2$, by
\begin{eqnarray*}
{\rm Cov}(T_k(\textbf{u}),T_k(\textbf{v}))&=&\sum_{j,l\neq k}[(\lambda_k - \lambda_j)(\lambda_k - \lambda_l)]^{-1}\mathbb{E}\left[\left< Z, \varphi_k \otimes \varphi_j \right>\left< Z, \varphi_k \otimes \varphi_l \right>\right] \varphi_j(\textbf{u})\varphi_l(\textbf{v})\\
&=& \sum_{j,l\neq k}[(\lambda_k - \lambda_j)(\lambda_k - \lambda_l)]^{-1}\mathbb{E}\left[\left< Z\otimes Z, (\varphi_k \otimes \varphi_j)\otimes (\varphi_k \otimes \varphi_l)\right>\right] \varphi_j(\textbf{u})\varphi_l(\textbf{v})\\
&=&\sum_{j,l\neq k}[(\lambda_k - \lambda_j)(\lambda_k - \lambda_l)]^{-1}\left[\left< \mathbb{E}(Z\otimes Z), (\varphi_k \otimes \varphi_j)\otimes (\varphi_k \otimes \varphi_l)\right>\right] \varphi_j(\textbf{u})\varphi_l(\textbf{v}).
\end{eqnarray*}
Since  $\Lambda= \mathbb{E}(Z\otimes Z)$, we have
\begin{equation}\label{cov}
{\rm Cov}(T_k(\textbf{u}),T_k(\textbf{v}))=\sum_{j\neq k}\sum_{l\neq k}\frac{\left< \Lambda, {\varphi}_k\otimes{\varphi}_j\otimes{\varphi}_k\otimes{\varphi}_j\right>}{(\lambda_k -\lambda_j)(\lambda_k -\lambda_l)}\varphi_j(\textbf{u})\varphi_l(\textbf{v}),
\end{equation}

Observe that the processes $T_{k,n}, k=1,\ldots, K$ are independent. This follows from Lemma 4.1 and the fact that the principal component functions $\varphi_k$ are independent.  Thus, the process $\mathbb{Z}_{xn}$ is a finte linear combination of independent processes, each of them converging in distribution to a Gaussian limit. %$\sum_{k=1}^K\alpha_k(x)\sqrt{n}[\hat{\varphi}_k(u,v) -\varphi_k(u,v) ]$
Then, $\mathbb{Z}_{xn}$ converges in distribution to the random element $\mathbb{Z}_x(u,v)=\textstyle\sum_{k=1}^K\alpha_k(x)T_k $, which is also a Gaussian process with covariance function equal to the sum of the covariances of the processes $T_k$ given in \eqref{cov}, where each term is multiplied by the factor $\alpha^2_k(x)$.\\

Now, we state in the following theorem  the weak convergence of our  empirical conditional copula process \(\{\sqrt{n}[\hat{C}_x(u, v) - C_x(u, v)]: 0 \leq u, v \leq 1\}\) to a Gaussian limit process.\\

\noindent\textbf{Theorem 4.1.}
Assume that the conditions of Propositions 4.1 and 4.2 hold. If, moreover, \(\sqrt{n}\tau_n \longrightarrow 0\) as \(n \longrightarrow \infty\), then we have
$$\sqrt{n}[\hat{C}_x(u,v) -C_x(u,v)]\stackrel{d}{\longrightarrow}\mathbb{W}_x(u,v),$$
where \(\mathbb{W}_x(u,v):=\mathbb{G}(u,v)+\mathbb{Z}_x(u,v)\) is a Gaussian process with covariance function
\begin{eqnarray*}
\Sigma[(u,v),(u',v')] &:=& \textup{Cov}(\mathbb{G}(u,v),\mathbb{G}(u',v')) + \textup{Cov}(\mathbb{Z}_x(u,v), \mathbb{Z}_x(u',v'))\\
&=& \bar{C}(u \wedge u', v \wedge v') - \bar{C}(u,v) \bar{C}(u',v')\\
&+& \sum_{k=1}^K \alpha_k^2(x) \sum_{j \neq k} \sum_{l \neq k} \frac{\left< \Lambda, \varphi_k \otimes \varphi_j \otimes \varphi_k \otimes \varphi_j \right>}{(\lambda_k - \lambda_j)(\lambda_k - \lambda_l)} \varphi_j(u,v) \varphi_l(u',v').
\end{eqnarray*}
\\

\textbf{Proof.} Let's recall the decomposition \eqref{dec2}:
$$\sqrt{n}[\hat{C}_x(u,v) -C_x(u,v)]=\mathbb{G}_n(u,v) + \mathbb{Z}_{xn}(u,v)+ O_P(\sqrt{n}\tau_n).$$
Now, we have to show that the sum process $\mathbb{G}_n+\mathbb{Z}_{xn}$ converges to a Gaussian process. Since both $\mathbb{G}_n$ and $\mathbb{Z}_{xn}$ converge in distribution to Gaussian processes $\mathbb{G}$ and $\mathbb{Z}_{x}$ respectively, it suffices to establish their asymptotic independence. For this purpose, we observe that the process $\mathbb{G}_n(u,v)$, given by \eqref{EMP}, depends on the variables $\varepsilon_{1i}$ and $\varepsilon_{2i}$, which are by definition independent of $X_i$ (see Lemma 2.1). On the other hand, we can express the process $\mathbb{Z}_{xn}(u,v)$  by using  observations $C_{X_i}$, whose randomness only depends on the $X_i$'s. That  is for every $(u,v)\in [0,1]^2$,
\begin{eqnarray*}
\mathbb{Z}_{xn}(u,v)&=&\sum_{k=1}^K\alpha_k(x)T_{k,n}(u,v)\\
&=&\sum_{k=1}^K\alpha_k(x)\sum_{j:j\neq k}{(\lambda_k -\lambda_j)^{-1}}\left< Z_n,{\varphi}_k\otimes{\varphi}_j\right>{\varphi}_j(u,v)\\
&=&\sqrt{n}\sum_{k=1}^K\sum_{j:j\neq k}\alpha_k(x){(\lambda_k -\lambda_j)^{-1}}\left< \hat{\Gamma}-\Gamma,{\varphi}_k\otimes{\varphi}_j\right>{\varphi}_j(u,v)\\
&=&\frac{1}{\sqrt{n}}\sum_{i=1}^n\sum_{k=1}^K\sum_{j:j\neq k}\alpha_k(x){(\lambda_k -\lambda_j)^{-1}}\\
&&\quad\times\left< (C_{X_i}-\tilde{\bar{C}}_n)\otimes(C_{X_i}-\tilde{\bar{C}}_n)-\Gamma,{\varphi}_k\otimes{\varphi}_j\right>{\varphi}_j(u,v)\\
&:=&\frac{1}{\sqrt{n}}\sum_{i=1}^n\phi(C_{X_i},u,v).
\end{eqnarray*}
Since $\varepsilon_{1i}$ and $\varepsilon_{2i}$ are independent of $X_i$, then they are independent of $C_{X_i}$ as well. Therefore, $\mathbb{G}_n(u,v)$ is independent of $\mathbb{Z}_{xn}(u,v)$ for all $(u,v)\in[0,1]^2$. Consequently, the sum process $\mathbb{G}_n +\mathbb{Z}_{xn}$ converges in distribution to $\mathbb{G} + \mathbb{Z}_x$, which is a Gaussian process with covariance function:
\begin{eqnarray*}
\Sigma[(u,v),(u',v')]&:=&\textup{Cov}(\mathbb{G}(u,v),\mathbb{G}(u',v'))+\textup{Cov}(\mathbb{Z}_x(u,v),\mathbb{Z}_x(u',v'))\\
 &=&\bar{C}(u\wedge u', v\wedge v')- \bar{C}(u,v) \bar{C}(u',v')\\
 &+&\sum_{k=1}^K\alpha_k^2(x)\sum_{j\neq k}\sum_{l\neq k}\frac{\left< \Lambda, {\varphi}_k\otimes{\varphi}_j\otimes{\varphi}_k\otimes{\varphi}_j\right>}{(\lambda_k -\lambda_j)(\lambda_k -\lambda_l)}\varphi_j(u,v)\varphi_l(u',v')
\end{eqnarray*} 
for all $(u,v)$ and $(u',v')$ in $[0,1]^2$. \qquad\qquad$\Box$

\newpage
\begin{center}
\section*{Appendix}
\end{center}
\textbf{Lemma A.1.}
For all $j=1,2,\dots,$ $$\left<\hat{\varphi}_j-\varphi_j,\varphi_j\right>=-\frac{1}{2}\|\hat{\varphi}_j-\varphi_j\|^2.$$
For all $j=1,2,\dots,$ $k=1,2,\dots,$ and $k\neq j,$ $$\left<\hat{\varphi}_j-\varphi_j,\varphi_k\right>=\frac{n^{-1/2}}{\hat{\lambda}_j-\lambda_k}\left<Z_n,\hat{\varphi}_j\otimes\varphi_k\right>,$$ for all $\hat{\lambda}_j\neq \lambda_k$.\\

\textbf{Proof.}
The definition of the norm in $L^2$, we have
$$\|\hat{\varphi}_j-\varphi_j\|^2= \left<\hat{\varphi}_j-\varphi_j,\hat{\varphi}_j-\varphi_j\right>.$$ 
By adding with $2\left<\hat{\varphi}_j-\varphi_j,\varphi_j\right>$, 
\begin{eqnarray*}
2\left<\hat{\varphi}_j-\varphi_j,\varphi_j\right>+\|\hat{\varphi}_j-\varphi_j\|^2&=&\left<\hat{\varphi}_j-\varphi_j,2\varphi_j\right>+\left<\hat{\varphi}_j-\varphi_j,\hat{\varphi}_j-\varphi_j\right>\\
&=&\left<\hat{\varphi}_j-\varphi_j,\varphi_j+\hat{\varphi}_j\right>\\
&=&\left<\hat{\varphi}_j,\varphi_j\right>-\left<\varphi_j,\varphi_j\right>+\left<\hat{\varphi}_j,\hat{\varphi}_j\right>-\left<\varphi_j,\hat{\varphi}_j\right>.
\end{eqnarray*}
Since $\left<\varphi_j,\varphi_j\right>=\|\varphi_j \|=1$, $\left<\hat{\varphi}_j,\hat{\varphi}_j\right>=\|\hat{\varphi}_j \|=1$ and that the inner product is symmetric, we obtain
$$2\left<\hat{\varphi}_j-\varphi_j,\varphi_j\right>+\|\hat{\varphi}_j-\varphi_j\|^2=0.$$
Through simple calculations, we can show that: 
$$\hat{\lambda}_k(\hat{\varphi}_j-\varphi_j) =\Gamma (\hat{\varphi}_j-\varphi_j)+\hat{\Gamma}(\hat{\varphi}_j)-\Gamma (\varphi_j)-(\hat{\lambda}_j-\lambda_j)\varphi_j.$$
Thus,
\begin{eqnarray*}
\hat{\lambda}_j\left<\hat{\varphi}_j-\varphi_j,\varphi_k\right>&=&\left<\Gamma,(\hat{\varphi}_j-\varphi_j)\otimes\varphi_k\right>+\left<\hat{\Gamma},\hat{\varphi}_j\otimes \varphi_k\right>-\left<\Gamma,\hat{\varphi}_j\otimes \varphi_k\right>\\
&=&\lambda_k\left<\hat{\varphi}_j-\varphi_j,\varphi_k\right>+\left<\hat{\Gamma}-\Gamma,\hat{\varphi}_j\otimes\varphi_k\right>.
\end{eqnarray*}
Which implies that
$$\left<\hat{\varphi}_j-\varphi_j,\varphi_j\right>=(\hat{\lambda}_j-\lambda_j)^{-1}\left<\hat{\Gamma}-\Gamma,\hat{\varphi}_j\otimes\varphi_j\right>=\frac{n^{1/2}}{\hat{\lambda}_j-\lambda_j}\left<Z_n,\hat{\varphi}_j\otimes\varphi_j\right>.\qquad\Box$$

%\newpage
%\bibliographystyle{plainnat}
%\bibliography{Biblio}
%\nocite{craiu2009parametric}
%\nocite{acar2011dependence}
%\nocite{abegaz2012semiparametric}
%\nocite{djaloud2024nonparametric}
%\nocite{li2019novel}
%\nocite{sabeti2013bayesian}

\end{document}